\documentclass[10pt]{amsart}\usepackage{amssymb} \usepackage{amscd}

\usepackage{amssymb}
\usepackage{amscd}
\usepackage{epsfig}
\usepackage[all]{xy}

\newcommand{\Z}{\mathbb Z}

\newcommand{\hm}{\mbox{Hom}}

\newcommand{\To}{\longrightarrow}

\newtheorem{thm}{Theorem}[section]
\newtheorem{lem}[thm]{Lemma}
\newtheorem{cor}[thm]{Corollary}

{\theoremstyle{definition} \newtheorem{rem}[thm]{Remark}
} {\theoremstyle{definition}
\newtheorem{defn}[thm]{Definition}}

{\theoremstyle{remark} }

\title{Weakly Exact Categories and the Snake lemma}
\author{Amir Jafari}
\begin{document}
\begin{abstract}
We generalize the notion of an exact category and introduce weakly exact categories. A proof of the snake lemma in this general setting is given. Some applications are given to illustrate how one can do homological algebra in a weakly exact category.
\end{abstract}
\maketitle
\tableofcontents
\section{Introduction}
Short and long exact sequences are among the most fundamental concepts in mathematics. The definition of an exact sequence is based on kernel and cokernel, which can be defined for any category with a zero object $0$, i.e. an object that is both initial and final. In such a category, kernel and cokernel are equalizer and coequalizer of:
\[\begin{xy}\xymatrix{
A\ar@<0.5ex>[r]^f\ar@<-0.5ex>[r]_0& B}\end{xy}.\]
Kernel and cokernel do not need to exist in general, but by definition if they exist, they are unique up to a unique isomorphism.
As is the case for any equalizer and coequalizer, kernel is monomorphism and cokernel is epimorphism. In an abelian category the converse is also true: any monomorphism is a kernel and any epimorphism is a cokernel. In fact any morphism $A\to B$ factors as $A\to I\to B$ with $A\to I$ a cokernel and $I\to B$ a kernel. This factorization property for a category with a zero object is so important that together with the existence of finite products, imply that the category is abelian (see 1.597 of \cite{FS}).

 A short exact sequence in an abelian category is a sequence of morphisms:
$$A\stackrel{i}\To B\stackrel{p}\To C$$
such that $i$ is a kernel of $p$ and $p$ is a cokernel of $i$. This definition makes sense for any category with a zero object and we call such sequences special. However if we let all special sequences to be short exact, the nice properties of exact sequences fail to extend unless we either restrict this definition or restrict the categories by putting extra conditions.

In the first approach we take a class of special sequences that satisfy certain carefully chosen axioms and call them short exact sequences. Such an approach for an additive category was considered by Quillen in \cite{Q}, and it has its predecessors e.g. in Heller \cite{H}. In the second approach
we assume that the category satisfies certain carefully chosen axioms (that are modeled on the categories of groups, rings, Lie algebras, etc.) and we take all special sequences to be short exact. Such an approach leads to the notions like Barr-exact, sequential, homological or semi-abelian categories , see \cite{BB} and the references therein.

Once we have the notion of a short exact sequence, one can easily define a long exact sequence.
\begin{defn} A long sequence
$$\dots\To A_{i-1}\To A_{i}\To A_{i+1}\To\dots$$
is exact if every morphism  $A_{i-1}\to A_i$ factors as $A_{i-1}\to Z_i\to A_i$ such that $Z_i\to A_i\to Z_{i+1}$ are short exact for all $i$. If the sequence is finite from one or both ends we also assume that the first $A_i\to Z_{i+1}$ sits in a short exact sequence $Z_i\to A_i\to Z_{i+1}$, and similarly the last $Z_i\to A_{i}$ sits in a short exact sequence $Z_i\to A_i\to Z_{i+1}$. \end{defn}

Among the consequences of exactness, an important result that is of constant use in homological algebra is the snake lemma:
\begin{lem}
In an abelian category, for the commutative diagram below, with exact rows:
\[\begin{xy}\xymatrix{
&A_1\ar[d]_{f_1}\ar[r]^{\phi_1}& A_2\ar[d]_{f_2}\ar[r]^{\phi_2}& A_3\ar[d]^{f_3}\ar[r]& 0\\
0\ar[r]& B_1\ar[r]_{\phi_1'}& B_2\ar[r]_{\phi_2'}& B_3}
\end{xy}\]
there is a natural\footnote{Natural here has a precise meaning: Given another diagram as above with $A_i'$, $B_i'$ and $f_i'$, and morphisms between them, making all squares commutative, we should have a commutative diagram:
\[\begin{xy}\xymatrix{
ker(f_3)\ar[r]^{\delta}\ar[d]& coker(f_1)\ar[d]\\
ker(f_3')\ar[r]^{\delta'}& coker(f_1')}
\end{xy}\]
}
connecting morphism $\delta: ker(f_3)\to coker(f_1)$ that makes the sequence
$$ker(f_1)\to ker(f_2)\to ker(f_3)\stackrel{\delta}\to coker(f_1)\to coker(f_2)\to coker(f_3)$$
exact. The maps between the kernels and the cokernels arise by their universal property. Furthermore if the top row is short exact then $0\to ker(f_1)\to ker(f_2)$ is exact and if the bottom row is short exact then
$coker(f_2)\to coker(f_3)\to 0$ is exact.
\end{lem}
The common proofs of this lemma use a method called chasing the elements. This can be applied for concrete categories, like abelian groups or more generally modules over a ring. The general case can be deduced using the embedding theorem that asserts any small abelian category can be embedded in the category of modules over a ring. There are several natural proofs without the use of elements in the literature, see e.g. \cite{FHH}.

In this note, we consider certain classes of special sequences in a general category with a zero object, that allow us to prove this lemma in a more general setting.
\section{Weakly Exact Category}
In this section we define a weakly exact category. These are categories with a zero object and certain class of morphisms, called deflation. deflations have kernel and are cokernel of their kernel. Hence if $p: B\to C$ is a deflation and $i:A\to B$ is its kernel then:
$$A\stackrel{i}\To B\stackrel{p}\To C$$
is a special sequence. We call such sequences short exact. A kernel of a deflation is called an inflation. Note that a deflation is epimorphism and an inflation is monomorphism. The terminology of admissible epimorphism and admissible monomorphism is also used for deflation and inflation. As in the introduction, we can also define long exact sequences in a weakly exact category.

The following axioms are very similar to Heller's axioms for an exact additive category, see appendix B of \cite{B} or \cite{H}. The difference is that we do not assume inflations are closed under composition (and hence the definition is not self dual), also axiom 3 is weaker than its counter part in Heller's axioms. Of course the major difference is that we no longer assume that our category is additive. The author does not know of any non-trivial example of a self dual weakly exact category that is not additive. However the following definition makes the category of groups, rings, Lie algebras, etc into a weakly exact category, if we take all (regular) epimorphisms to be deflation.
\begin{defn}
A weakly exact category is a category with a zero object and a collection of morphisms called deflations satisfying the following axioms.
\begin{enumerate}
\addtocounter{enumi}{-1}
\item Any isomorphism and any morphism $A\to 0$ is a deflation.
\item Any deflation has a kernel and is cokernel of its kernel. The kernel of a deflation is called an inflation.
\item The class of deflations is closed under composition.
\item For two composable morphism $f$ and $g$ if $gf$ and $f$ are deflation then $g$ is a deflation.
\item (3 by 3 lemma)  In the following commutative diagram with short exact columns if the second and the third rows are short exact then the first row is short exact.
        \[\begin{xy}\xymatrix{
        A_1\ar[d]\ar[r]&A_2\ar[r]\ar[d]&A_3\ar[d]\\
        B_1\ar[d]\ar[r]&B_2\ar[r]\ar[d]&B_3\ar[d]\\
        C_1\ar[r]&C_2\ar[r]&C_3\\
        }
        \end{xy}\]
\end{enumerate}
\end{defn}
\begin{lem} Axiom 4 is self dual, i.e. in the diagram above with short exact rows, if the first and second columns are short exact then the third column is short exact.
\end{lem}
{\bf Proof.} The morphism $B_3\to C_3$ is a deflation by axiom 2 and 3. So let $A_3'\to B_3$ be its kernel. By properties of kernel we have a map from $A_2\to A_3'$ which makes the diagram above with $A_3$ replaced with $A_3'$ commutative. Its columns and the second and the third rows are short exact, hence axiom 3 implies that $A_1\to A_2\to A_3'$ is short exact. Since both $A_3$ and $A_3'$ are cokernel of $A_1\to A_2$, there is an isomorphism $A_3\to A_3'$ such that the composite $A_3\to A_3'\to B_3$ is a kernel for $B_3\to C_3$, hence $A_3\to B_3\to C_3$ is short exact.\qed
\section{The Snake Lemma}
In this section we will prove:
\begin{thm} (The snake lemma) In a weakly exact category, given a commutative diagram
\begin{equation}\begin{xy}\xymatrix{
A_1\ar[d]_{f_1}\ar[r]^{\phi_1}& A_2\ar[d]_{f_2}\ar[r]^{\phi_2}& A_3\ar[d]^{f_3}\\
 B_1\ar[r]_{\phi_1'}& B_2\ar[r]_{\phi_2'}& B_3}
\end{xy}\end{equation}
where the rows are short exact and $f_i$'s are admissible, i.e they have a factorization
$$A_i\stackrel{f_i'} \To I_i\stackrel{f_i''}\To B_i$$
with $f_i'$ a deflation with kernel $k_i: K_i\to A_i$ and $f_i''$ an inflation with cokernel $c_i:B_i\to C_i$, there is a natural connecting morphism $\delta: K_3\to C_1$ such that $$K_1\stackrel{\psi_1}\To K_2\stackrel{\psi_2}\To K_3\stackrel{\delta}\To C_1\stackrel{\psi_1'}\To C_2\stackrel{\psi_2'}\To C_3$$ is exact and $\psi_1$ is an inflation and $\psi_2'$ is a deflation. Here $\psi_i$ and $\psi_i'$ are the natural morphisms induced by $\phi_i$ and $\phi_i'$.
\end{thm}
{\bf Proof.} This proof is inspired by the proof of the snake lemma for a sequential category given in \cite{Bo}.
First observe that by axioms 2 and 3, $\psi_2'$ is a deflation. By the universal properties of kernel and cokernel there is a unique morphism $\phi_2'': I_2\to I_3$, such that the following diagram is commutative:
\[\begin{xy}\xymatrix{
A_2\ar[r]^{\phi_2}\ar[d]_{f_2'}&A_3\ar[d]^{f_3'}\\
I_2\ar@{.>}[r]^{\phi_2''}\ar[d]_{f_2''}&I_3\ar[d]^{f_3''}\\
B_2\ar[r]^{\phi_2'}&B_3}\end{xy}\]

By axioms 2 and 3, $\phi_2''$ is a deflation and let $\phi_1'': J_1\to I_2$ be its kernel. Then by axiom 4, the sequence $J_1\to B_1\to ker(\psi_2')$ in the following diagram is short exact.
\begin{equation}\begin{xy}\xymatrix{
        J_1\ar[d]_{\phi_1''}\ar[r]^{\epsilon}&B_1\ar[r]^{c_2'}\ar[d]^{\phi_1'}&ker(\psi_2')\ar[d]^i\\
        I_2\ar[d]_{\phi_2''}\ar[r]^{f_2''}& B_2\ar[r]^{c_2}\ar[d]^{\phi_2'}& C_2\ar[d]^{\psi_2'}\\
        I_3\ar[r]^{f_3''}&B_3\ar[r]^{c_3}&C_3\\
        }
        \end{xy}\end{equation}
Note that since the composite $C_1\stackrel{\psi_1'}\To C_2\stackrel{\psi_2'}\To C_3$ is trivial $\psi_1'$ factors through $ker(\psi_2')$ via a morphism $\alpha: C_1\to ker(\psi_2')$. Since $\alpha c_2=c_2'$ and both $c_2$ and $c_2'$ are deflations, it follows that $\alpha$ is a deflation. This shows that $C_1\to C_2\to C_3$ is exact.
\\
Consider the following diagram
  \begin{equation}\begin{xy}\xymatrix{
        I_1\ar@{=}[d]\ar[r]^{\epsilon'}&T\ar[r]^j\ar[d]^{\epsilon}&ker(\psi_1')\ar[d]^{i'}\\
        I_1\ar[d]\ar[r]^{f_1''}&B_1\ar[r]^{c_1}\ar[d]&C_1\ar[d]^{\psi_1'}\\
        0\ar[r]&ker(\psi_2')\ar@{=}[r]&ker(\psi_2')\\
        }
        \end{xy}\end{equation}

where $T$ is the kernel of $\psi'_1c_1=c_2\phi_1'$.  Axiom 4 again implies that the first row is a short exact sequence. Since in diagram (2), $c_2\phi_1'=ic_2'$, we can identify $T\to B_1$ with $J_1\to B_1$.
Let $\pi_1: P\to A_2$ be the kernel of $f_3'\phi_2: A_2\to I_3$ or equivalently the kernel of $\phi_2''f_2': A_2\to I_3$. It fits in the following commutative diagrams:
\[\begin{xy}\xymatrix{
    A_1\ar@{=}[d]\ar@{.>}[r]^{\theta_1}&P\ar@{.>}[r]^{\theta_2}\ar[d]^{\pi_1}&K_3\ar[d]^{i}&&K_2\ar@{.>}[r]^{k_2'}\ar@{=}[d]&P\ar@{.>}[r]^{p}\ar[d]^{\pi_1}&J_1\ar[d]^{\phi_1''}\\
    A_1\ar[d]\ar[r]_{\phi_1}&A_2\ar[r]_{\phi_2}\ar[d]& A_3\ar[d]^{f_3'}&&K_1\ar[d]\ar[r]_{k_2}&A_2\ar[r]_{f_2'}\ar[d]&I_2\ar[d]^{\phi_2''}\\
    0\ar[r]&I_3\ar@{=}[r]&I_3&&0\ar[r]&I_3\ar@{=}[r]&I_3}
    \end{xy}
    \]
By axiom 4 the first rows are short exact. Consider the commutative diagram:
\[\begin{xy}\xymatrix{
A_1\ar[r]^{\theta_1}\ar[d]^{f_1'}& P\ar[r]^{\theta_2}\ar[d]^p & K_3\ar@{.>}[d]^{\delta'}\\
I_1\ar[r]_{\epsilon'}&J_1\ar[r]_j& ker(\psi_1')}\end{xy}\]
To show that $p\theta_1=\epsilon'f_1'$ it is enough to show $\epsilon p\theta_1=\epsilon\epsilon'f_1'$. Now an easy chase in the diagrams show that both sides are $f_1$.
We define $\delta: K_3\stackrel{\delta'}\To ker(\psi_1')\stackrel{i'}\hookrightarrow C_1$. Note that since $p$, $j$ and $\theta_2$ are deflations, axiom 2 and 3 imply that $\delta'$ is a deflation. This shows that $K_3\to C_1\to C_2$ is exact. Finally if we apply axiom 4 to the diagram above, we deduce that
$$K_1\to K_2\to ker(\delta')$$
is a short exact sequence. The theorem is proved.\qed
\begin{lem} Let $f:A\to B$ and $g:B\to C$ be two morphisms. If $gf$ and $g$ are inflation then $f$ is an inflation.\end{lem}
{\bf Proof.} Let $B\stackrel{g}\To C\stackrel{p}\To D$ and $A\stackrel{gf}\To C\stackrel{q}\To D'$ be short exact sequences. Consider the following diagram:
\[\begin{xy}\xymatrix{
A\ar[r]^{gf}\ar[d]&C\ar[r]^q\ar[d]^p&D'\ar@{.>}[d]^{p'}\\
0\ar[r]&D\ar@{=}[r]&D}\end{xy}\]
The left square commutes since $fg$ factors through the kernel of $p$, i.e. $B$. The morphism $p'$ exists by universal property of cokernel. Now axiom 3 implies that $p'$ is a deflation. So we can apply axiom 4 and deduce that $A\stackrel{f}\to B\to ker(p')$ is a short exact sequence so $f$ is an inflation.\qed

\begin{rem} If in the statement of the snake lemma we compose $A_1\to A_2$ with a deflation $A_1'\to A_1$ and compose $B_2\to B_3$ with an inflation $B_3\to B_3'$, we get natural morphisms $K_1'\to K_1$ and $C_3\to C_3'$ where $K_1'$ is the kernel of $A_1'\to A_1\to B_1$ and $C_3'$ is the cokernel of $A_3\to B_3\to B_3'$. It follows easily from the axioms that $K_1'\to K_1$ is a deflation and $C_3\to C_3'$ is an inflation (for this we need to use the previous lemma). Hence we have a long exact sequence $K_1'\to K_2\to K_3\to C_1\to C_2\to C_3'$.
\end{rem}
\section{The $3\times 3$ Lemma}
In this section we show that axiom 4 follows from the following two axioms:
\\
(4a) The pullback of a deflation along a deflation exists and is a deflation.
\\
(4b) Given a commutative diagram:
\[\begin{xy}\xymatrix{
A_1\ar[d]^{f_1}\ar[r]^{\phi_1}&A_2\ar[dr]^{\phi_2}\ar[d]^{f_2}\\
B_1\ar[r]_{\phi_1'}&B_2\ar[r]_{\phi_2'}&A_3}\end{xy}\]
where $\phi_2$ and $\phi_2'$ are deflation and $\phi_1$ and $\phi_1'$ are their kernels, $f_1$ is a deflation if and only if $f_2$ is a deflation.

 These are usually easier to check in a category. However axiom 4 does not in general imply (4a) and (4b).

Axiom 4, says in the commutative diagram:
\begin{equation}\begin{xy}\xymatrix{
A_1\ar[d]_{f_1}\ar[r]^{\phi_1}& A_2\ar[d]_{f_2}\ar[r]^{\phi_2}& A_3\ar[d]^{f_3}\\
 B_1\ar[r]_{\phi_1'}& B_2\ar[r]_{\phi_2'}& B_3}
\end{xy}\end{equation}
where the rows are short exact and $f_i$'s are deflations with kernel $K_i$, the induced sequence $K_1\to K_2\to K_3$  is short exact.

Consider the pullback of $\phi_2'$ along $f_3$:
\[\begin{xy}\xymatrix{
     &P\ar[d]_{\pi_2}\ar[r]^{\phi_2''}& A_3\ar[d]^{f_3}\\
     B_1\ar[r]_{\phi_1'}\ar@{.>}[ur]^{\phi_1''}&B_2\ar[r]_{\phi_2'}&B_3}
     \end{xy}\]
Since $\phi_2'f_2=f_3\phi_2$, by the universal property of pullback, there is a unique morphism $\pi_1: A_2\to P$ such that $\pi_2\pi_1=f_2$ and $\phi_2''\pi_1=\phi_2$. Also since the kernel of a morphism and its pullback are isomorphic, we have a lift of the kernel $\phi_1'$ of $\phi_2'$ to $\phi_1'':B_1\to P$ which is a kernel of $\phi_2''$. We therefore get the following commutative diagram:
\[\begin{xy}\xymatrix{
     A_1\ar[r]^{\phi_1}\ar[d]^{f_1}&A_2\ar[dr]^{\phi_2}\ar[d]^{\pi_1}\\
     B_1\ar[r]^{\phi''_1}\ar[dr]_{\phi_1'}&P\ar[d]^{\pi_2}\ar[r]^{\phi_2''}& A_3\ar[d]^{f_3}\\
     &B_2\ar[r]^{\phi_2'}&B_3}
     \end{xy}\]
 Now by axiom 4b, $\pi_1$ is a deflation whose kernel is  $k_1':K_1\stackrel{k_1}\to A_1\stackrel{\phi_1}\to A_2$. Similarly the kernel of $f_3$, $k_3: K_3\to A_3$ has a unique lift $k_3':K_3\to P$ as a kernel of $\pi_2$. Note that $\pi_2$ is a deflation by axiom 4a. Finally by axiom (4b) the commutative diagram:
\[\begin{xy}\xymatrix{
K_2\ar[r]^{k_2}\ar[d]^{\psi_2}&A_2\ar[d]^{\pi_1}\ar[dr]^{f_2}\\
K_3\ar[r]_{k_3'}&P\ar[r]_{\pi_2}&B_2}\end{xy}\]
implies that $\psi_2$ is a deflation whose kernel is the same as $\pi_1$, i.e. $K_1\stackrel{\psi_1}\to K_2$ hence the sequence $K_1\to K_2\to K_3$ is a short exact sequence.
For completeness we prove the full version of the 3 by 3 lemma:
\begin{lem} In a weakly exact category satisfying (4a) and (4b), given a commutative diagram:
\[\begin{xy}\xymatrix{
        A_1\ar[d]_{f_1}\ar[r]^{\phi_1}&A_2\ar[r]^{\phi_2}\ar[d]^{f_2}&A_3\ar[d]^{f_3}\\
        B_1\ar[d]_{g_1}\ar[r]^{\phi_1'}&B_2\ar[r]^{\phi_2'}\ar[d]^{g_2}&B_3\ar[d]^{g_3}\\
        C_1\ar[r]_{\phi_1''}&C_2\ar[r]_{\phi_2''}&C_3\\
        }
        \end{xy}\]
where all the columns and the first and the last rows are short exact and $\phi_2'\phi_1'=0$, then the middle row is short exact.
\end{lem}
{\bf Proof.} Consider the following pullback:
\[\begin{xy}\xymatrix{
P\ar[r]^{\theta_2}\ar[d]_{\psi_2}&C_2\ar[d]^{\phi_2''}\\
B_3\ar[r]_{g_3}&C_3}\end{xy}\]
As in the previous proof we get  morphisms $\theta_1:A_3\to P$ as kernel of $\theta_2$ and $\psi_1: B_2\to P$ such that $\psi_2\psi_1=\phi_2'$ such that the following diagram is commutative:
\[\begin{xy}\xymatrix{
A_2\ar[r]^{f_2}\ar[d]_{\phi_2}&B_2\ar[d]_{\psi_1}\ar[dr]^{g_2}\\
A_3\ar[r]^{\theta_1}\ar[dr]_{f_3}&P\ar[r]^{\theta_2}\ar[d]_{\psi_2}&C_2\ar[d]^{\phi_2''}\\
&B_3\ar[r]_{g_3}&C_3}\end{xy}\]

Since $\phi_2''$ and $\phi_2$ are deflations, by axiom (4a) and (4b) $\psi_2$ and $\psi_1$ are deflations and hence $\phi_2'=\psi_2\psi_1$ is a deflation. Let $\tilde{\phi}_1': B_1'\to B_2$ be its kernel. By axiom (4) the induced sequence $A_1\stackrel{f_1'}\To B_1'\stackrel{g_1'}\To C_1$ is short exact. Since $\phi_2'\phi_1'=0$ the morphism $\phi_1'$ factors through $\tilde{\phi}_1'$ via $\psi:B_1\to B_1'$. We have a commutative diagram:
\[\begin{xy}\xymatrix{
&B_1\ar[d]^{\psi}\ar[dr]^{g_1}\\
A_1\ar[r]^{f_1'}\ar[ur]^{f_1}& B_1'\ar[r]^{g_1'}&C_1}\end{xy}\]

Now by axiom (4b), $\psi$ is a deflation and since
\[\begin{xy}\xymatrix{
A_1\ar[r]^{f_1}\ar@{=}[d]&B_1\ar[d]^{\psi}\\
A_1\ar[r]_{f_1'}&B_1'}\end{xy}\]
is a pullback, kernel of $\psi$ is the same as the kernel of the identity of $A_1$, and hence is zero. So $0\to B_1\stackrel{\psi}\to B_2'$ is short exact i.e. $\psi$ is an isomorphism. Therefor $B_1\to B_2$ is a kernel of $B_2\to B_3$ and the middle row is short exact. \qed
\section{Examples}

Recall that a morphism that is coequalizer of two parallel arrows is called a regular epimorphism. For example any cokernel is a regular epimorphism. In a weakly exact category, any deflation being a cokernel is regular epimorphism. We make the following definition:
\begin{defn} A category with a zero object is said to be weakly homological, if the class of all regular epimorphisms satisfy the axioms for a weakly exact category. Obviously this is a the maximal class.
\end{defn}

For example any abelian category is weakly homological. In this case every epimorphism is regular.

Recall that a regular category is a category such that
\begin{enumerate}
\item All finite limits exists.
\item The class of regular epimorphisms is closed under arbitrary pullbacks.
\item For any morphism $f:X\to Y$ the kernel pair:
\[\begin{xy}\xymatrix{
X\times_Y X\ar@<0.5ex>[r]^{p_1}\ar@<-0.5ex>[r]_{p_0}& X}\end{xy}.\]
has a coequalizer $g:X\to Z$ and therefore $f=hg$ for a unique morphism $h:Z\to Y$.
\end{enumerate}
It is shown that $h$ is a monomorphism. So any morphism can be factored as a regular epimorphism followed by a monomorphism, in a unique way.
\begin{lem} A regular category with zero object is weakly homological if and only if every regular morphism is cokernel of its kernel.\end{lem}
{\bf Proof.} The "only if" part is clear, since axiom 1 should hold for all regular epimorphisms. To prove the other direction, first note that in a regular category axioms 0,2 and 3 hold (see appendix A.5 of \cite{BB}). Actually for axiom 3, it is enough to know that $gf$ is a regular epimorphism.
Also axiom 4a and half of 4b (if $f_2$ is regular epic then $f_1$ is) are true, since regular epics are stable under pullbacks. It remains to show the other half of 4b. Factor $f_2$ as $A_2\stackrel{f_2'}\To P \stackrel{f_2''}\To B_2 $, with $f_2'$ regular epic and $f_2''$ monic.
since $f_1$ and $f_2'$ are cokernel of their isomorphic kernels, by universal properties of cokernel we get a morphism $\phi_1'': B_1\to P$ making the following diagram commutative:

\[\begin{xy}\xymatrix{
A_1\ar[d]^{f_1}\ar[r]^{\phi_1}&A_2\ar[dr]^{\phi_2}\ar[d]^{f_2'}\\
B_1\ar[r]_{\phi_1''}\ar@{=}[d]&P\ar[r]_{\phi_2''}\ar[d]^{f_2''}&A_3\ar@{=}[d]\\
B_1\ar[r]_{\phi_1'}&B_2\ar[r]_{\phi_2'}&A_3
}\end{xy}\]
Since $f_2''$ is monic, $\phi_1''$ is kernel of $\phi_2''$ and therefore the lower right square is a pullback. Hence $f_2''$ is a deflation. And therefore $f_2=f_2''f_2'$ is a deflation.\qed

We give an example of a category, where every regular epimorphism is cokernel of its kernel.
\begin{lem} In any additive category with finite limits, any regular epic is cokernel of its kernel.\end{lem}
{\bf Proof.} Let $f:X\to Y$ be a regular epimorphism with kernel $k:K\to X$ and kernel pair $\begin{xy}\xymatrix{
X\times_Y X\ar@<0.5ex>[r]^{p_1}\ar@<-0.5ex>[r]_{p_0}& X}\end{xy}.$ Since $fp_0=fp_1$ hence $f(p_0-p_1)=0$ and hence $p_0-p_1$ factors through $k$. Therefore for $g:X\to Z$ if $gk=0$ then $g(p_0-p_1)=0$ and hence $g$ coequalizes $p_0$ an $p_1$. Since $f$ is the universal coequalizer of $p_0$ and $p_1$, $g$ factors uniquely through $f$, this is the definition of a cokernel for $k$. \qed

More generally, it is shown in \cite{Bo}, that in any protomodular category with a zero object, any regular epic is cokernel of its kernel and hence a regular protomodular  category with a zero object is weakly homological. Such a category is called homological in \cite{BB}. A category is Protomodular if for any pullback diagram:
\[\begin{xy}\xymatrix{
A'\ar[r]^{\phi'}\ar[d]_{f}&B'\ar[d]^{g}\\
A\ar[r]_{\phi}&B}\end{xy}\]
where $f$ and $g$ are split epic, i.e. we have $s:A\to A'$ and $t:B\to B'$ such that $fs=id_{A}$ and $gt=id_{B}$, the pair $(\phi',t)$ is jointly strongly epic. Only the following special case of this definition is needed to show any regular epic is cokernel of its kernel. Given a regular epic $f:X\to Y$ with kernel $k:K\to X$ and kernel pair $\begin{xy}\xymatrix{
X\times_Y X\ar@<0.5ex>[r]^{p_1}\ar@<-0.5ex>[r]_{p_0}& X}\end{xy},$ we get the following pullback diagram:
\[\begin{xy}\xymatrix{
K\times K\ar@<-0.5ex>[d]_{q_0}\ar@<0.5ex>[d]^{q_1}\ar[r]^{k'}& X\times_Y X\ar@<-0.5ex>[d]_{p_0}\ar@<0.5ex>[d]^{p_1}\\
K\ar[r]^k& X}\end{xy}\]
of split epic morphisms, where  $s:K\to K\times K$ and $t:X\to X\times_Y X$ are given by the diagonals. Then the pair $(k',t)$ is jointly epic.

\section{Applications}

    In this section we discuss several applications of the snake lemma, similar to those in homological algebra on abelian categories. But instead we work in a weakly exact category. We will consider chain complexes in such a category, i.e. sequences $\dots \to A_{i-1}\to A_i\to A_{i+1}\to\dots $ where the composition of successive arrows are zero. We have a natural definition for morphisms between such chain complexes. The morphisms $A_i\to A_{i+1}$ are called differentials and if they are admissible morphisms then the chain complex is called an admissible complex.
    \begin{defn}
    The cohomology of a chain complex $A\stackrel{f}\to B\stackrel{g}\to C$ with inflation $f$ and deflation $g$ is defined by one the following two ways:
    \begin{enumerate}
    \item Let $A\stackrel{f}\to B\stackrel{g'}\to C'$ be a short exact sequence, by properties of cokernel $g$ factors through $C'$:  $B\stackrel{g'}\to C'\to C$. By axiom 3, $C'\to C$ is a deflation and its kernel $H$ is the cohomology.
    \item Let $A'\stackrel{f'}\to B\stackrel{g}\to C$ be a short exact sequence, by properties of kernel $f$ factors through $A'$ : $A\to A'\stackrel{f'}\to B$. By lemma 3.2., $A\to A'$ is an inflation and its cokernel $H'$ is the cohomology.
        \end{enumerate}
        \end{defn}
        \begin{lem}
        With the notation of the previous definition, there is a natural morphism $H\to H'$ that is both deflation and inflation and hence an isomorphism. Furthermore the complex $A\to B\to C$ is exact if and only if $H=0$.
        \end{lem}
{\bf Proof.} Applying the snake lemma to the following diagram proves the first part:
\[\begin{xy}\xymatrix{
A\ar[r]\ar[d]&B\ar@{=}[d]\ar[r]&C'\ar[d]\\
A'\ar[r]&B\ar[r]&C\\
}
\end{xy}
\] The second part is easy and is left to the reader.\qed
\begin{defn} A differential object in $\mathcal C$ is an object $A$ together with an admissible morphism $d_A:A\to A$ such that $d_Ad_A=0$. A morphism between two differential objects $(A,d_A)$ and $(B, d_B)$ is a morphism $f:A\to B$ such that $fd_A=d_Bf$. We have therefore a category formed by differential objects denoted by $d\mathcal C$. The cohomology defined in the previous definition associates to a differential object $(A,d_A)$ the cohomology of the admissible chain $$Im(d_A)\hookrightarrow A\stackrel{d_A}\twoheadrightarrow Im(d_A).$$ It is easy to see that it defines a functor $H:d{\mathcal C}\to {\mathcal C}$.
\end{defn}

\begin{thm} Let $A \to A' \to A''$ is a short exact sequence of differential objects, then there is a natural morphism $H(A'')\to H(A)$ that makes the following triangle exact at each vertex:
\[\begin{xy}\xymatrix{
H(A)\ar[rr]&&H(A')\ar[dl]\\
&H(A'')\ar[ul]}
\end{xy}\]
\end{thm}
{\bf Proof.} We can decompose the admissible morphism $d:A\to A$ as $A\stackrel{d}\to I\hookrightarrow A$ with $A\to I$ a deflation and $I\to A$ an inflation. So we have short exact sequences $I\to A\to C$ and $K\to A\stackrel{d}\to I$. By properties of kernel and cokernel, the morphism $A\stackrel{d}\to I$ factors through $C$ : $A\to C\to I$ with $C\to I$ a deflation and similarly $I\to A$ factors through $K$ : $I\to K\to A$ with $I\to K$ an inflation, hence we have an admissible morphism $C\to K$ as composite of $C\to I\to K$. According to the previous lemma its kernel and cokernel are both $H(A)$. Such constructions can be done for $A'$ and $A''$. So we have the following commutative diagram:
\[\begin{xy}\xymatrix{
C\ar[r]\ar[d] & C'\ar[r]\ar[d]& C''\ar[d]\\
K\ar[r] & K'\ar[r] & K''}
\end{xy}\]
By snakes lemma applied to the diagram:
\[\begin{xy}\xymatrix{
A\ar[r]\ar[d]^d & A'\ar[r]\ar[d]^{d'}& A''\ar[d]^{d''}\\
A\ar[r] & A'\ar[r] & A''}
\end{xy}\]
we realize that $C\to C'\to C''\to 0$ is exact and $0\to K\to K'\to K''$ is exact. So if we apply the snake lemma (modified as in remark 3.3.) to the first diagram, we get the desired morphism $H(A'')\to H(A)$ making the triangle in the statement of the theorem exact at each vertex.\qed
\begin{cor} If $A\to B\to C$ is a short exact sequence of differential objects and if two of these differential objects are exact (i.e. their cohomology are zero) then so is the third one.
\end{cor}
\begin{defn} The $i$th cohomology of an admissible chain complex
$$\dots\To A_{k-1}\stackrel{d_{k-1}}\To A_k\stackrel{d_k}\To A_{k+1}\To\dots$$
Is either the cokernel of the inflation $Z_{i-1}\hookrightarrow K_i$ or the kernel of the deflation $C_{i-1}\stackrel{d_i}\to Z_i$ where $Z_i=Im(d_i)$, $K_i=ker(d_i)$ and $C_i=coker(d_i)$. It is denoted by $H^i(A_{\bullet})$.
\end{defn}
\begin{thm}
Let $A_{\bullet}\to A_{\bullet}'\to A_{\bullet}''$ be a sequence of chain morphisms between admissible chain complexes such that for each $i$, $A_i\to A_i'\to A_i''$ is short exact (we call these special sequences, point wise short exact). Then we have a long exact sequence of cohomologies:
$$\dots\to H^i(A_{\bullet})\to H^i(A_{\bullet}')\to H^i(A_{\bullet}'')\to H^{i+1}(A_{\bullet})\to H^{i+1}(A_{\bullet}')\to H^{i+1}(A_{\bullet}'')\to\dots$$
\end{thm}
{\bf Proof.} We have morphisms $C_{i-1}\to K_{i+1}$ as composite $C_{i-1}\to Z_i\to K_{i+1}$. The kernel of this morphism is $H^i$ and the cokernel is $H^{i+1}$. The result follows if we apply the snake lemma to the following diagram:
\[\begin{xy}\xymatrix{
C_{i-1}\ar[r]\ar[d] & C'_{i-1}\ar[r]\ar[d]& C''_{i-1}\ar[d]\\
K_{i+1}\ar[r] & K'_{i+1}\ar[r] & K''_{i+1}}
\end{xy}\]
The fact that $C_{i-1}\to C'_{i-1}\to C_{i-1}''\to 0$ and $0\to K_{i+1}\to K_{i+1}'\to K_{i+1}''$ are exact follows from the snake lemma as well.\qed
\begin{lem} The category of admissible chain complexes with the class of deflations being point wise deflation is a weakly exact category.\end{lem}  {\bf Proof.} The only nontrivial fact to be proved is to show that if $f:A_{\bullet}\to B_{\bullet}$ is point wise deflation, then its kernel $K_{\bullet}$ is an admissible chain complex. Let $A_i\to I_i\to A_{i+1}$ and $B_i\to I'_i\to B_{i+1}$ be admissible decomposition of the differentials, then by properties of kernel and cokernel and axioms 2 and 3 we have deflation morphisms $I_i\to I_i'$. Let $I_i''\to I_i$ be its kernel. Then the morphisms $K_i\to K_{i+1}$ factor through $I_i''$ and by axiom 3 it's easy to show that the decomposition $K_i\to I_i''\to K_{i+1}$ is an admissible decomposition for the differentials of $K_{\bullet}$.\qed
\section{Additive Weakly Exact Categories}
In this section we want to extend the notion of chain homotopy . We need the following definitions:
\begin{defn} A super structure on a weakly exact category $\mathcal A$, is a choice of a functor $\epsilon:{\mathcal A}\to {\mathcal A}$, such that:
\begin{enumerate}
\item $\epsilon^2=1$.
\item $\epsilon$ is identity on objects.
\item If $A\stackrel{f}\To B\stackrel{g}\To C$ is short exact then $A\stackrel{\epsilon(f)}\To B\stackrel{\epsilon(g)}\To C$ is short exact.
\end{enumerate}
\end{defn}
\begin{defn} Let $\mathcal A$ be an additive weakly exact category with an (additive) super structure  $\epsilon$. A differential object compatible with the super structure in $\mathcal A$ is a pair of an object $A$ and a morphism $d_A:A\to A$, such that $d_A^2=0$ and $\epsilon(d_A)=-d_A$. It is
admissible if $d_A$ is an admissible morphism. For admissible differential objects $H(A)$ is defined in the previous section.
\end{defn}
With these definitions we can prove the following:
\begin{lem} If $(A,d_A)$ and $(B,d_B)$ are two differential objects compatible with the super structure in a category as above. Then the morphism:
$$d:\hm(A,B)\To\hm(A,B)$$
$$df=fd_A-d_B\epsilon(f)$$
is a differential, i.e. $d^2=0$. Let $Z(A,B)=ker(d)$ , ${\Bbb B}(A,B)=Im(d)$. If $A$ and $B$ are admissible, then we have a natural morphism
$$Z(A,B)/{\Bbb B}(A,B) \To \hm(H(A),H(B))$$
If we let the objects to be differential objects compatible with the super structure and morphisms to be $Z(A,B)/{\Bbb B}(A,B)$, then we get an additive category called the homotopy category of super differential objects and is denoted by $H(d{\mathcal A})$. We let the full subcategory of admissible differential objects to be $H(d^a{\mathcal A})$. The morphism defined above, gives us a functor
$$H:H(d^a{\mathcal A})\To {\mathcal A}.$$
\end{lem}
{\bf Proof.} Note that
$$d^2f=(fd_A-d_B\epsilon(f))d_A-d_B(-\epsilon(f)d_A+d_Bf)=0.$$
Now if $f\in ker(d)$ then $fd_A=d_B\epsilon(f)$. So $f:A\to B$ maps the image of $d_A$ to the image of $d_B$.  Applying $\epsilon$ we get $\epsilon(f)d_A=d_Bf$, so $f$ maps the kernel of $d_A$ to the kernel of $d_B$. Furthermore if $f=dg=gd_A-d_B\epsilon(g)$ then $f$ maps kernel of $d_A$ to the image of $d_B$ and this implies that the induced map $H(A)\to H(B)$ is zero. Therefore we get the desired map $Z(A,B)/{\Bbb B}(A,B)\To \hm(H(A),H(B))$. The second part is easy, we need to check that if $f\in Z(A,B)$ and $g\in {\Bbb B}(B,C)$ then $gf\in {\Bbb B}(A,C)$ and similarly if $h\in Z(C,D)$ then $hg\in {\Bbb B}(B,D)$.\qed
\\
An important notion for the homotopy category of differential objects is that of a quasi isomorphism. We start with the admissible case:
\begin{defn} A morphism $f:A\to B$ in $H(d^a{\mathcal A})$ of admissible differential objects is said to be a quasi isomorphism if the induced map $H(A)\to H(B)$ is an isomorphism.\end{defn}
This definition can not be applied for the general case, since $H(A)$ and $H(B)$ might not be defined. However we have the following lemma:
\begin{lem} If $f:A\to B$ is a morphism between two differential objects in $Z(A,B)$ then for any differential object  $X$, the map:
\[\hm(X,A)\To\hm(X,B)\]
given by sending $g$ to $fg$ is a morphism of differential abelian groups. So it induces a map $H(\hm(X,A))\to H(\hm(X,B))$. If $f\in {\Bbb B}(A,B)$ then this map is zero and we therefore get:
$$Z(A,B)/{\Bbb B}(A,B)\To \hm(H(\hm(X,A)),H(\hm(X,B)))$$
\end{lem}
\begin{defn} A morphism $f:A\to B$ in $H(d{\mathcal A})$ is weakly quasi isomorphism if the induced map $H(\hm(X,A))\to H(\hm(X,B))$ is an isomorphism for all differential objects $X$.
\end{defn}
\begin{lem} If $f:A\to B$ in $H(d^a{\mathcal A})$ is a quasi-isomorphism then it is weakly quasi isomorphism in $H(d{\mathcal A})$. \end{lem}

We would like to apply these ideas to the category of chain complexes. So let us begin with a weakly exact additive category $\mathcal A$.
\begin{defn} The category of graded objects in $\mathcal A$, is a category whose objects are sequences $(A_i)_{i\in\Z}$ of objects of $\mathcal A$, indexed by integers and a morphism $f$ between $(A_i)$ and $(B_i)$ is a sequence of morphisms $f_i:A_i\to B_{i+k}$ for a fixed $k$, which is called the degree of $f$ and is denoted by $\deg(f)$. We use ${\mathcal A}^{\bullet}$ to denote this category.\end{defn}
It is easy to see that ${\mathcal A}^{\bullet}$ is a weakly exact category if we let a deflation $f$ be a point-wise deflation, i.e. all $f_i:A_i\to B_{i+k}$ are deflations. We also have a super structure $\epsilon$:
$$\epsilon(f)=(-1)^{\deg(f)}f.$$
An (admissible) chain complex, as in the previous section, is simply a differential object in ${\mathcal A}^{\bullet}$ with an (admissible) differential of degree 1. Although ${\mathcal A}^{\bullet}$ is not additive (we can add two morphisms, only when they have same degrees\footnote{This is just a technical issue, instead of taking the union of morphisms of different degrees, we could have taken their coproduct in the category of abelian groups.}), the constructions of the previous lemma goes through so long as $d_A$ and $d_B$ have the same degree. If $\deg(d_A)=\deg(d_B)=1$, $\hm(A,B)$ will be a chain complex of abelian groups. So we let $H^0(d{\mathcal A})$ be a category whose objects are differential objects in ${\mathcal A}^{\bullet}$ with a differential of degree 1 and the hom-set between $(A^{\bullet},d_A)$ and $(B^{\bullet},d_B)$ is $H^0(\hm(A,B))$. If we assume the differential is admissible we get the subcategory $H^0(d^a{\mathcal A}^{\bullet})$. This is the homotopy category of (admissible) chain complexes in $\mathcal A$. It is an additive category. We also use the notation ${\mathcal K}({\mathcal A})$ (resp ${\mathcal K}^a({\mathcal A})$) for this category. As in the lemma we have a functor:
\[H: {\mathcal K}^a({\mathcal A})\To {\mathcal A}^{\bullet}\]
The $n$th component of this functor is denoted by $H^n$.
\\
As before we have a definition for (weakly) quasi isomorphism in ${\mathcal K}^a({\mathcal A})$ and ${\mathcal K}({\mathcal A})$.
\\

One should have a derived categories ${\mathcal D}^a({\mathcal A})$ and ${\mathcal D}({\mathcal A})$, with a functor ${\mathcal K}^{(a)}({\mathcal A})\To {\mathcal D}^{(a)}({\mathcal A})$ such that any (weakly) quasi isomorphism is sent to an isomorphism and it is universal, i.e. for any other functor ${\mathcal K}^{(a)}({\mathcal A})\To {\mathcal C}$ with the above property, it should factor through the derived category. There is still another way to define quasi isomorphism, a morphism $f\in Z^0(A,B)$ is a quasi-isomorphic if its cone is an exact sequence. We leave the study of these matters for future and end our note here.

\end{document}